\documentclass[12pt]
{article}
\usepackage{amssymb,amsmath, mathtools}

\usepackage[latin1]{inputenc}
\usepackage[T1]{fontenc}
\usepackage{play}

\usepackage{array}   
\newcolumntype{L}{>{$}l<{$}} 

\usepackage{caption}

\usepackage{graphicx}
\DeclareGraphicsExtensions{.pdf,.png,.jpg}

\graphicspath{ {images/} }

\usepackage{lipsum}
 \usepackage{enumitem}

\usepackage{amsmath}
\usepackage{amsthm}
\usepackage{amsfonts}
\usepackage{bbm}
\usepackage{amssymb}
\usepackage{graphicx}
\usepackage{color}

\usepackage{wasysym, marvosym}
\usepackage{lmodern}
\usepackage{subfig}

\usepackage{fancybox}
\newlength{\querylen}
\newcommand{\prob}{\mathbb{P}}

\newcommand{\jmp}[1]{\\[#1pt]}
\newcommand{\jump}{\\[-5pt]}
\newcommand{\sbygone}{{\tt bygone}}
\newcommand{\snext}{{\tt next}}
\newcommand{\strap}{{\tt trap}}

\newcommand{\successprofile}{{\boldsymbol p}}
\newcommand{\priordist}{{\boldsymbol \pi}}
\newcommand{\nth}[2]{$#1^{\text{#2}}$}

\usepackage[flushleft]{threeparttable}

\newtheorem{thm}{Theorem}
\newtheorem{lemma}[thm]{Lemma}

\newtheorem{cor}[thm]{Corollary}

\newtheorem{assertion}[thm]{Proposition}
\theoremstyle{definition}

\theoremstyle{remark}

\begin{document}
\title{Trapping  the Ultimate Success} 
\author{Alexander Gnedin ~~~and~~~Zakaria Derbazi\\ {\it\small Queen Mary, University of London}}
\maketitle

\begin{abstract}

\noindent{}We introduce a betting game, where the gambler aims to guess the last success epoch from past observed data.
The player may bet on the event that no further successes occur, or choose a {\it{}trap} which is any span of future times. Winning is achieved if the last success turns out to be the only one falling in  
the trap. The game is closely related to the sequential decision problem of   maximising the probability of stopping on the last success in a finite sequence of trials.
We use this connection to analyse the problem of stopping at the last record for trials paced by a P{\'o}lya-Lundberg process with log-series distribution of the total  number of trials.

\end{abstract}

\section{Introduction}

Suppose a series of inhomogeneous Bernoulli  trials with given {\it  profile} of success probabilities $\successprofile =(p_k,~k\geq1)$
 is paced randomly in time by some independent  point process.
 As the outcomes and epochs
of the first $k\geq0$ trials get  known at some time $t$,
 the gambler is asked to bet on the time of the last success. 
The gambler is allowed to choose from three strategies: 1) a $\strap$ strategy where the gambler wins when the last success epoch gets isolated in a proper set of future times, that is, 
it falls in the trap while no other success epoch occurs. 2) A $\snext$ strategy, where winning is achieved if exactly one success happens in the future and 3) a $\sbygone$ strategy, where the gambler wins if no further successes occur.

A classic profile is related to the random records model where the trials can be uniquely ranked and  the exchangeability of ranks entails $p_k=1/k$. For such profile, trapping is an instance of a stopping strategy in the best choice problem where the objective is to recognise the overall best trial, {\it the last record}, at the moment it occurs \cite{BG, Browne, BR, CZ, AnoK, Stewart, TamakiWang}.
Other  choices  of $\successprofile$ are suggested by random combinatorial structures and many other areas where inhomogeneous Bernoulli trials play an eminent  role \cite{PB}.

The analogous  trapping game with discrete time is amenable to study  by means of the optimal stopping theory for Markov chains. As a consequence, the state space for the sequence of successes is just the set of  natural numbers. Thus, every Markovian stopping time coincides with a trapping strategy determined by a set of integers. 
The problem with fixed number of trials and general
$\boldsymbol p$ has been discussed in previous research \cite{odds, Ribas}, and \cite{PresmanSonin} treats the best choice problem with a random number of trials.
However, these previous setups are different from the continuous time game. Both the index of trial and its time are important decision variables.


Regarding the pacing point process, we shall assume that it is mixed binomial without multiple points. The assumption entails that the pair $(t,k)$ is a sufficient statistic summarising the observed data before time $t$.
The setting covers the wide class of mixed Poisson processes and many others.
In a nutshell, the pacing process is characterised by the {\it prior} distribution  $\priordist$ of the total number of trials, and some background continuous distribution of i.i.d. `arrivals'.
Without loss of generality,  the model is standardised by assuming that the arrivals are uniformly distributed. That is to say, whenever the number of trials is $n$,  they are paced at locations of the uniform order statistics on $[0,1]$.

The most obvious instance of a trapping strategy is a {\it $z$-strategy}. 
For $0<z<1$, this leaves the $(1-z)$ {\it proportion} of the remaining time to trap the last success.  The edge value $z=0$ corresponds to the action {\tt next}. 
We will give a simple condition on $\successprofile$ to ensure that a  $z$-strategy is optimal among all trapping strategies.

A question of  central interest in this paper is the characterisation of pairs ($\successprofile, \priordist$) which admit that
 for some states $(t,k)$, the action
$\sbygone$ outperforms $\snext$ but a proper trapping is better still.
The question is motivated by the optimal stopping problems, in which a gambler's online strategy is an arbitrary adapted stopping time, and  the objective is to stop at the last success epoch. 
If the stopping problem belongs to the so-called {\it monotone case} \cite{CRS}, the optimal strategy is {\it myopic},  that is prescribing to 
stop  at the earliest success epoch when $\sbygone$ becomes more beneficial than $\snext$. 
Therefore, in the monotone case trapping cannot be better than  
both $\sbygone$ and $\snext$. However, trapping can be used  to assess if the stopping problem belongs to the monotone case.

It is inherent in the model to consider each  prior  within the context of  a family of  power series distributions 
\begin{eqnarray}\label{psd}
\pi_n=c(q)w_n q^n\,~,n\geq0,
\end{eqnarray}
with given shape weights $(w_n)$, and  scale
 parameter $q>0$. This allows one to define the critical cutoffs for trapping and 
 the stopping problem in terms of roots of certain power series in the variable $x=(1-t)q$.

The random records model with geometric prior has a special feature that the point process of record epochs is Poisson.
Then the optimal strategy has a single cutoff approaching $1/e$ as $q\to1$ \cite{BrussRogers, BrussSam1, BrussSam2}.
The limit form, commonly  called the  $1/e${\it -strategy}, coincides with 
$z$-strategy  for $z=1/e$  referred to the decision time  $t=0$.
It is known that the $1/e$-strategy stops at the last success with probability at least $1/e$, provided the number of trials is non-zero, and this bound
is the best possible \cite{Bruss84, BrussSam1, G1/e}.
Recently,   it was observed  \cite{BR} that the $1/e$-strategy is not optimal for the problem with trials occurring at times of a linear birth process. 
Here, we will cast the model of \cite{BR} in the context of the  log-series prior, apply $z$-strategies, and show  that
the stopping problem does not belong to the monotone case.

\section{Definitions} 

\subsection{The probability model}

Let $\priordist$ be a power series distribution (\ref{psd}) with weights
$$ w_0\geq 0, ~~~w_n>0{\rm ~~~for~~}n\geq 1.$$
The associated mixed binomial process 
on the unit interval
 is an orderly  counting process $(N_t, ~t\in[0,1])$ with the uniform order statistic property.
The process, can also be seen as a time inhomogeneous pure-birth process, with  transition rate  expressible through the generating function of $(w_n)$, see \cite{Puri}.
The posterior distribution of the number of trials yet to occur is again a power series distribution
\begin{equation}\label{poster}
\pi(j\,|t,k):=\prob(N_1-N_t=j|N_t=k)= f_k(x) {k+j\choose j}    w_{k+j} x^j,~~j\geq0,
\end{equation}
with scale variable
\begin{equation}\label{x-t}
x=(1-t)q
\end{equation}
and a normalisation function  $f_k(x)$. 
The conditioning relation (\ref{poster}) 
appears in many statistical problems related  to censored data. 

In principle, instead of considering a family of  processes $(N_t)$ with parameter $q$, we could deal with one Markov process defined as function of the `size' variable (\ref{x-t}),
where $q>0$ assumes values within the range of convergence of  $\sum_n w_n q^n$.
We prefer  not to adhere to this viewpoint, as  the `real time' parameter is more  intuitive. Nevertheless,  we will switch back and forth between $t$ and $x$, as $x$ is more suitable
 for power series work.

The arrivals are assumed to be uniformly distributed due to the nice self-similarity features. That is, 
conditionally on $N_t=k$
\begin{itemize}
\item[(i)] The point processes of trials on $[0,t)$ and $(t,1]$ are independent,
\item[(ii)]  $\big(N_{t+s/(1-t)}-N_t,~s\in[0,1]\big)$ is a mixed binomial process on $[0,1]$,  with the number of trials distributed according to  (\ref{poster}).
\end{itemize}

Let $\successprofile=(p_k, ~k\geq1)$ be a profile of success probabilities. 
We assume that 
$$0\leq p_1\leq 1,~~~0\leq p_k<1~ ~~{\rm for~}~k>1 \text{ and }\sum_{k=1}^\infty p_k=\infty.$$
The trial at {\it index} $k$ (\nth{k}{th} trial) is a success with probability $p_k$, independently of other trials  and the pacing process.
Thus, the point process of success epochs is obtained from $(N_t)$ by thinning out the  \nth{k}{th} point with probability $1-p_k$.
Typically, the point process of success epochs is not Poisson, nor even Markovian.

We denote $(t,k)$ the state of  the counting process, meaning the  event $N_t=k$, and write $(t,k)^\circ$ if the \nth{k}{th} trial occurs at time $t$ and it is a success.

\subsection{The trapping game and stopping problem}

A single round of the  trapping game played in the generic state $(t,k)$ is the following.
The player chooses either a proper subset of the interval $(t,1]$, 
$\snext$ or $\sbygone$.
 A $z$-strategy corresponds to the interval with endpoints $t+z (1-t)$ and $1$.
Such interval is called {\it final}, and the left endpoint is called cutoff.
The gambler's  objective is  to choose an admissible action to maximise the probability of isolating  the last success epoch  from other successes.

For the trapping game, it is irrelevant whether the state is $(t,k)$ or $(t,k)^\circ$.
The game in state $(t,k)$ can be reduced to the game in state $(0,0)$, by assigning to
(\ref{poster}) the role of  a prior, and truncating the profile of success probabilities. 
In state $(0,0)$, the trap of $z$-strategy is just the final interval $(z,1]$.

A reason to  consider the state $(t,k)$ as a variable, is the connection with the following optimal stopping problem (as mentioned in the Introduction):

Consider the increasing filtration of sigma-algebras induced by the natural information flow, 
so that  the data available at time $t$  comprises the location of trials on $[0,t]$ and their outcomes.
Let $\tau$ be an adapted stopping time, viewed as online strategy of the gambler.  
For a given succession of the trials,  
the range of $\tau$ is the set of success epochs or $1$. 
The gambler wins a pound if the last success epoch is $(\tau,N_\tau)^\circ$, otherwise there is no payoff. In particular, there is no payoff in the event $\tau=1$.
The objective is to maximise the winning probability (equal to the expected payoff).

Stopping time is said to be Markovian if  the decision  in state $(t,k)^\circ$ only depends on the state, but not on the trials before time $t$, or the initial state.
A trapping strategy can be seen as a randomised  stopping time  initiated in  some state $(t_0,k_0)$ or $(t_0,k_0)^\circ$.
It is non-Markovian because it depends on the initial state.
In particular,  for 
 $z$-strategy the stopping condition involves   
 $t\geq t_0+z(1-t_0)$.

\section{The fixed-$n$ game}

The trapping game with  fixed number of trials is not trivial itself.  
This can be seen as a game of informed gambler who learns  the number of future trials at time of the decision.

\subsection{Discrete time}
Suppose the gambler in state $(t,k)$ learns that there are $j$ trials yet to occur, making the total to $n=k+j$.
The number of successes in unseen trials $k+1,\cdots, n$
 has probability generating function 
$$\lambda\mapsto \prod_{m=k+1}^n (1-p_m+p_m\lambda)= \left(1+\lambda \sum_{i=k+1}^n \frac{p_i}{1-p_i} \right)   \prod_{m=k+1}^n (1-p_m)+O(\lambda^2).$$
The probability of zero successes is
$$s_0(k+1,n):= \prod_{m=k+1}^n (1-p_m),$$
and the probability of exactly one success is
$$s_1(k+1,n):= \sum_{i=k+1}^n \frac{p_i}{1-p_i}\prod_{m=k+1}^n (1-p_m)= s_0(k+1,n)\sum_{i=k+1}^n \frac{p_i}{1-p_i}.$$

There is an obvious recursion relationship between $s_0$ and $s_1$
$$s_1(k,n)=(1-p_k)s_1({k+1},n)+p_ks_0(k+1,n),$$
which we can write as
\begin{eqnarray}\nonumber
s_1(k,n)-s_1(k+1,n)&=&p_k \{s_0(k+1,n)-s_1(k+1,n)\}\\
&=& p_ks_0(k+1,n)\left(1-\sum_{i=k+1}^n \frac{p_i}{1-p_i} \right) \label{recur}
\end{eqnarray}
Since the sequence
\begin{eqnarray}\label{sum1}
1-\sum_{i=k+1}^n \frac{p_i}{1-p_i},~~~0\leq k\leq n-1,
\end{eqnarray}
has at most one variation of sign, namely its sign pattern is
$$-,\cdots, -,\geq 0,+,\cdots,+,$$
It follows that:
\begin{itemize}
\item[(i)] $s_1(\cdot,n)$ is unimodal with at most  two (adjacent) maximum locations,
\item[(ii)] The modes are non-decreasing in $n$.
\end{itemize}
For $n$ fixed, we also have that:
\begin{itemize}
\item[(iii)] The mode is precisely the minimal location where  {\tt bygone} starts outperforming {\tt next}.
\end{itemize}

Let $k^*$ be the mode of $s_1(\cdot \,,n)$ and $A^*=\{k^*,\cdots,n\}$. A well known fact of the optimal stopping theory \cite{odds} is that $A^*$ is optimal, in the sense that 
no other set $A\subset\{1,\cdots,n\}$ isolates the last success in $n$ trials with higher probability. The following is a direct variational proof:
Clearly, $n\in A$ is necessary  for $A$ to be optimal. By induction, 
suppose we have shown that $\{k+1,\cdots,n\}\subset A$.  Including $k$ impacts  the winning chance by 
$$c \,p_k  \{s_0(k+1,n)-s_1(k+1)\}$$
where $c\geq  0$ 
depends on $A\cap \{1,\cdots,k-1\}$ only. If this is non-negative then $k$ should be included, otherwise not. 
Hence $A=A^*$.

\paragraph{Monotonicity in $n$.} 
We show next that the sign of 
\begin{equation}\nonumber
\max_k s_1(k,n)-\max_k s_1(k,n+1) 
\end{equation}
is the same as the sign of
$p_{k^*+1}-p_{n+1}$. In particular, the winning probability is non-increasing  if $p_n\downarrow$.
To argue this point, we identify  $A^*=\{k^*+1,\cdots,n\}$  with a stopping strategy in discrete time.
By increasing the number of trials by one,  the mode may either remain the same or increase by one.

Firstly, compare how $A^*$ performs in the $n$-problem,
with the stopping set $B:=\{k^*+1,\cdots,n+1\}$ applied  in the problem with $n+1$ trials. 
Clearly, strategies $A^*$ and $B$ only differ if  
the \nth{(n+1)}{st} trial is a success, and the number of successes in trials $k^*+1,\cdots,n$ is either $1$ or $0$.
Thus, the winning probabilities differ by
$$(s_1(k^*+1,n)-s_0(k^*+1,n))p_{n+1}= \left(1-\sum_{i=k^*+1}^n \frac{p_i}{1-p_i}\right)s_0(k^*+1,n)\leq 0.$$

Secondly, compare $A^*$ with the other option, $C:=\{k^*+2,\cdots,n,n+1\}$. The difference of winning probabilities of $A^*$ in $n$-problem and $C$ in $(n+1)$-problem  has four component probabilities:
\begin{enumerate}
\item[(a)]  $p_{k^*+1} s_0(k^*+2,n)(1-p_{n+1})$
that \nth{(k^*+1)}{st} trial is a success, $A$ wins, $B$ loses,
\item[(b)] $(1-p_{k^*+1}s_1(k^*+2,n)p_{n+1}$ that \nth{(k^*+1)}{st}  trial is a failure,  $A$ wins, $B$ loses,
\item[(c)] $p_{k^*+1}s_1(k^*+2,n)(1-p_{n+1})$ that\nth{(k^*+1)}{st}  trial is a success,  $A$ loses, $B$ wins,
\item[(d)] $(1-p_{k^*+1})s_0(k^*+2,n)p_{n+1}$ that \nth{(k^*+1)}{st}  trial is a failure,  $A$ loses, $B$ wins.
\end{enumerate}
After simplification, (a)+(b)-(c)-(d) becomes
$$\left(1-\sum_{i=k^*+2}^n \frac{p_i}{1-p_i}\right)(p_{k^*+1}-p_{n+1}),$$
which has the same sign as $p_{k^*+1}-p_{n+1}$ because the first factor is non-negative by the optimality of $A^*$.

\subsection{Fixed $n$, trapping in continuous time}

In the elementary  continuous-time scenario,  a fixed number $n$  of trials occur  at uniformly sampled locations on $[0,1]$.
In state $(t,k)$,
the trapping probability for $z$-strategy is a Bernstein polynomial in $z$,
\begin{eqnarray}  \label{S1}
S_1(k, n; z):=\sum_{j=0}^{n-k} {n-k\choose j}    z^j (1-z)^{n-k-j}s_1(k+j+1, \,n).
\end{eqnarray}
Replacing $s_1$ by $s_0$ in this formula gives the probability denoted $S_0(k,n;z)$, that none of the successes gets trapped by the $z$-strategy,
with $S_0(k,n;0)$ equal to the probability to win with $\sbygone$.
The dependence on $t$ is void, since conditionally on $k$ arrivals before $t$, there is a non-random number $n-k$ of arrivals uniformly paced in $(t,1]$.

Note that $s_0(k+1,n)=S_0(k,n; 0)$ and $s_1(k+1,n) = S_1(k,n; 0)$.
The  form of the optimal stopping strategy in the fixed-$n$ discrete-time problem 
and
the theorem about excluding randomised stopping times \cite{CRS} imply that
\begin{equation}\label{trivial}
S_0(k,n; 0) \geq S_1(k,n; 0)\Longrightarrow  S_1(k,n; 0)=\max_z S_1(k, n; z).
\end{equation}
That is to say, trapping is ineffective if $\sbygone$ is better than $\snext$. This holds for $k, n$ satisfying
$$s_0(k+1,n)\geq s_1(k+1,n) \Longleftrightarrow  \sum_{j=k+1}^n \frac{p_j}{1-p_j}\leq 1.$$
Replacing a final interval by any other trap does not change the conclusion.

From the unimodality of $s_1(\cdot,n)$ and the shape-preserving properties of the Bernstein polynomials (see \cite{Bernstein}, Theorem 3.3), it follows that
(\ref{S1}) is unimodal.  Therefore, a unique strategy exists which is optimal among the $z$-strategies.
Mimicking the discrete-time variational argument, it will be shown in the following that  other traps (Borel sets) cannot be better.

\begin{thm}\label{T1} The optimal trapping strategy on $n$ trials is a $z$-strategy, where $z$ is  the unique mode of $S_1( k,n;\cdot)$.
The mode is $0$ in the case {\rm (\ref{trivial})}, and otherwise $z\in(0,1)$.
\end{thm}
\begin{proof} 
To ease notation,  we consider the state $(0,0)$, which is sufficient.
There is certainly a final interval that belongs to the optimal trap, because near the end of the time interval
the probability of two or more successes is negligible. Now, suppose $[z,1]$ belongs to the trap and we are assessing if the length element $[z-{\rm d}z,z]$ is worth including.
The change  of the winning probability  due to the inclusion is some factor depending on the structure of the trap within $[0,z-{\rm d}z]$ multiplied by the following

\begin{eqnarray}\label{incr}
 \sum_{k=1}^n {n-1\choose k-1} z^{k-1} (1-z)^{n-k}  p_k  \{s_0(k+1,n)-s_1(k+1,n)\}  \,n\,{\rm d}z&=&\\ \nonumber
 (1-z)^n \sum_{k=1}^n {n-1\choose k-1}  \left(\frac{z}{1-z}\right)^{k}   p_k  \{s_0(k+1,n)-s_1(k+1,n)\}  \,n\,{\rm d}z.
\end{eqnarray}
By (\ref{recur}),
in the variable $z/(1-z)$ the polynomial $\sum (\cdots)$ has at most one  variation of sign in the coefficients. Applying Descartes' rule of signs, we see that the polynomial has at most one positive root.
This  implies that the optimal trap is a final interval with cutoff coinciding with the root, 
or $[0,1]$ if there are no roots.\jump

It remains to check that the root, if any, coincides  with the mode of 
$$S_1(0,n; z)=\sum_{k=0}^{n} {n\choose k}    z^k    (1-z)^{n-k}s_1(k+1, \,n).$$
Indeed, we have for the derivative using (\ref{recur})

\begin{eqnarray*}
D_z S_1(0,n; z)&=& \sum_{k=1}^n  {n-1\choose k-1}n z^{k-1}(1-z)^{n-k}s_1(k+1,n)-
\sum_{k=0}^{n-1} {n-1\choose k}n z^k(1-z)^{n-k-1}s_1(k+1,n)\\
&=&\sum_{k=1}^n (\cdots)-\sum_{k=1}^n {n-1\choose k-1}n z^{k-1}(1-z)^{n-k}s_1(k,n)\\
&=&\sum_{k=1}^n  {n-1\choose k-1}n z^{k-1}(1-z)^{n-k}   \{s_1(k+1,n)-s_1(k,n)\}\\
&=&\sum_{k=1}^n  {n-1\choose k-1}n z^{k-1}(1-z)^{n-k}  p_k\{s_1(k+1,n)-s_0(k+1,n)\},
\end{eqnarray*}
Which is the negative of the polynomial in  (\ref{incr}). This gives the desired conclusion.
\end{proof}

\subsection{Examples}

\paragraph{The best choice problem.} 
In the random records model, 
the  formula $p_k=1/k$ for probability of record is a consequence of the exchangeability of ranks of the trials.
The Bernstein polynomials satisfy
$$S_1(k,n;z)\to -z\log z, ~~~n\to\infty,$$
and the convergence is uniform.   The sequence of modes converges to $1/e$.

In the case $k=0$, the Bernstein polynomial can be written in the form of a Taylor polynomial
$$S_1(0,n;z)=1-z- \sum_{j=2}^{n} \frac{(1-z)^j}{j(j-1)},$$
which decreases pointwise to  $z\mapsto -z\log z$ as $n$ increases. As observed in \cite{Bruss84}, the modes increase monotonically to $1/e$ and also
$$\max_z S_1(0,n;z) \downarrow 1/e.$$
These facts underline the minimax property of the stopping strategy with a single cutoff $1/e$,
known as the $1/e$-strategy
 \cite{Bruss84, Bruss88, BrussSam2}. See \cite{G1/e} for a recent  analysis of strategic dominance of this and other minimax strategies.

For $k>0$, the above nice monotonicity properties are no longer valid, the minimax value is below $1/e$ and  the $1/e$-cutoff  strategy is not minimax.
  This is seen already in the case $k=1$, where the Bernstein polynomials have   alternative representations
\begin{eqnarray}\nonumber
S_1(1, n;z)&=& \frac{n-1}{n}-\sum_{j=2}^{n-1} \frac{(n-j)(1-z)^j}{n j (j-1)}\\
    &=&S_1(0,n; z) +\sum_{j=1}^{n-1} \frac{(1-z)^{j+1}}{n\,j}  - \frac{(1-z)}{n}.
\label{k=1Bern}
\end{eqnarray}
The first formula is derived by conditioning on the highest rank $j$ of the trials that occur on $[0,z]$.

\begin{figure}
\centering
\subfloat{{\includegraphics[scale=0.22]{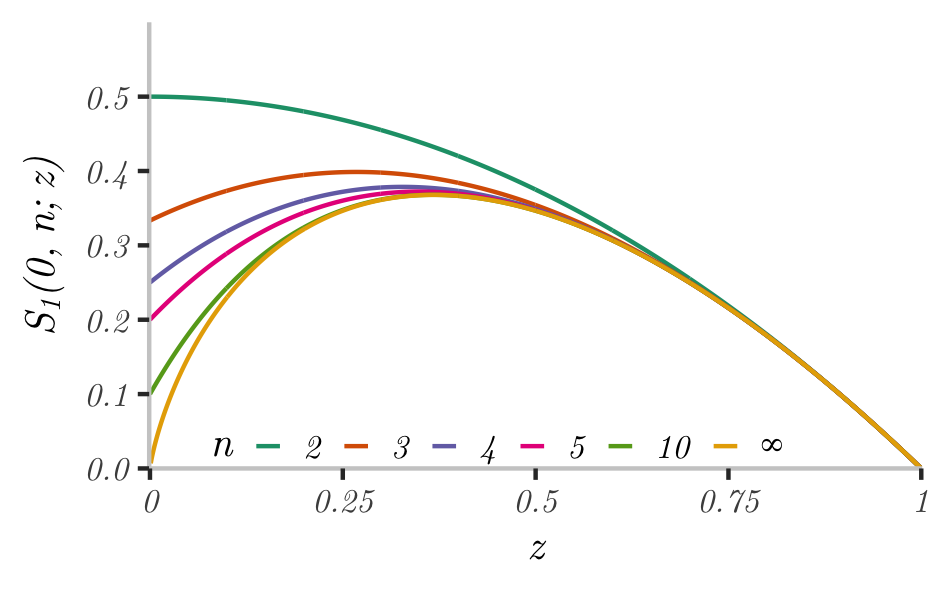} }}%
\,
\subfloat{{\includegraphics[scale=0.22]{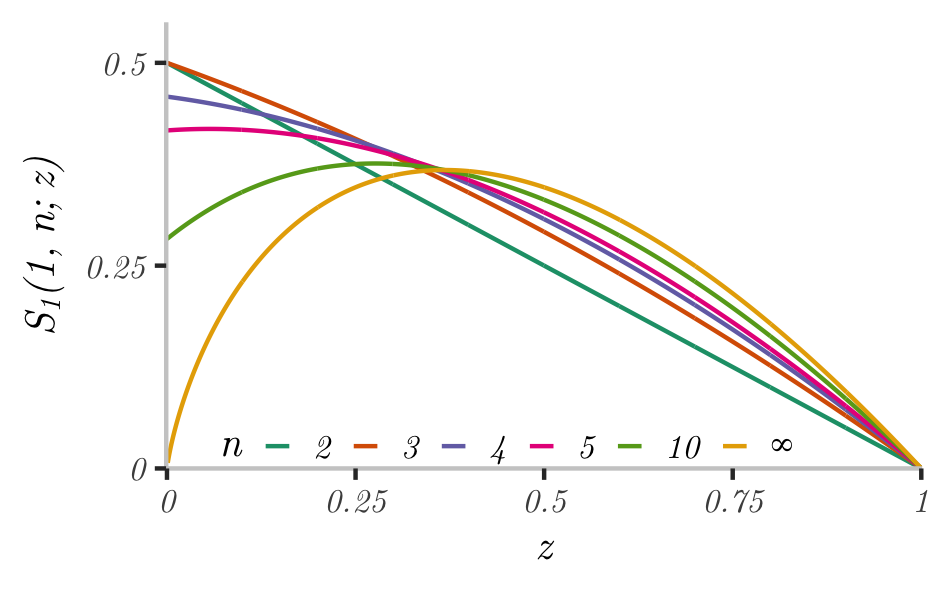} }}%
\caption{Bernstein polynomials for $p_k=1/k$.}%
\label{fig:bernstein}%
\end{figure}

\paragraph{The Karamata-Stirling profile.}  
The profile 
$$p_k=\frac{\theta}{\theta+k-1}, ~~k\geq 1,$$
with parameter $\theta>0$, plays a central role in the combinatorial structures related to the Ewens sampling formula for random partitions \cite{ABT}.
The term  {\it Karamata-Stirling law}  was coined in \cite{Bingham} for the distribution of the number of successes with these probabilities.
The number of successes in trials $k+1,\cdots,n$ has probability generating function
$$\lambda\mapsto \frac{(k+\theta\lambda)_{n-k}}{(k+\theta)_{n-k}}.$$
As $n\to\infty$,   $S_1(0,n;z)\to -\theta z^\theta\log z$ and the modes converge to $e^{-1/\theta}$. 
The shapes vary considerably with $\theta$. For 
 $\theta$ large, the minimax trapping value is close to zero.

\begin{figure}[h]
\centering
\subfloat{{\includegraphics[scale=0.22]{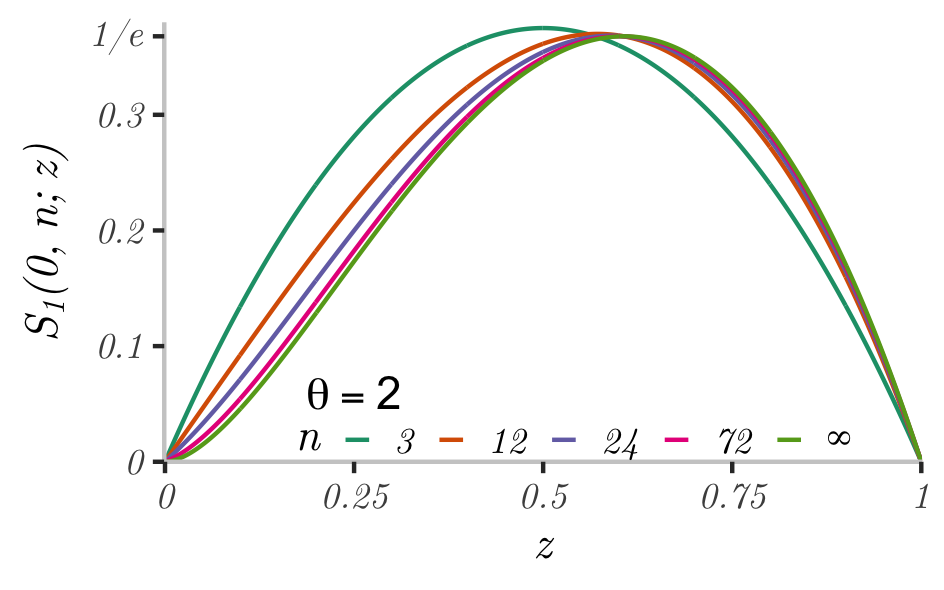} }}%
\,
\subfloat{{\includegraphics[scale=0.22]{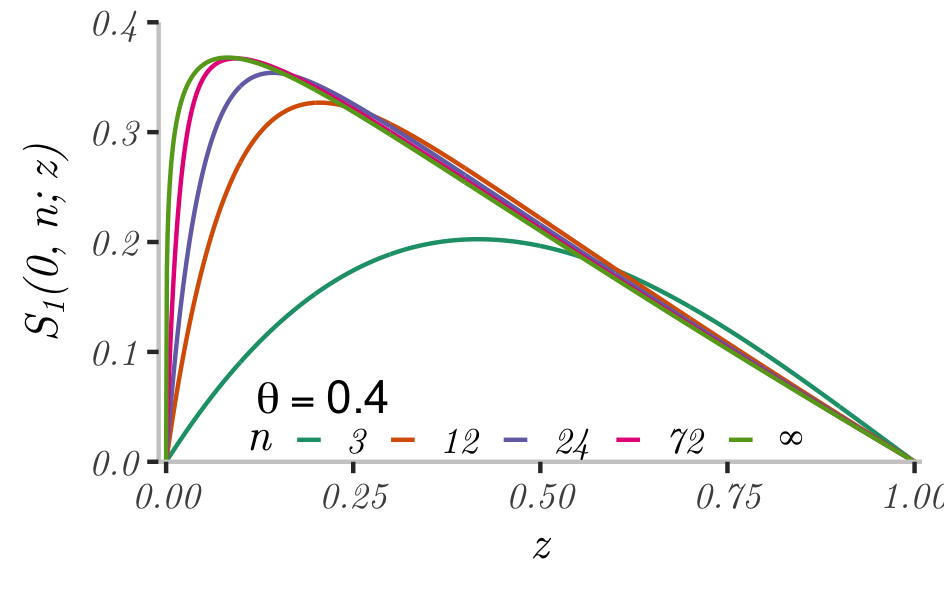} }}%
\caption{Bernstein polynomials for $p_k=\theta/(\theta+k-1)$.}
\label{fig:bernstein2}%
\end{figure}

\vspace*{-20pt}

\section{Random number of trials}

The best choice discrete-time problem with random number of trials was pioneered in \cite{PresmanSonin}.
The following features can be readily  extended to  the general  profiles $\successprofile$.
The sequence of success epochs  is  a  Markov chain on non-negative integers, and a stopping strategy can be identified with a set of integers.\jmp{3}
In general, the optimal stopping set $A^*$ is not a gap-free integer interval, it is rather comprised of `stopping islands' whose number and 
configuration depend on the prior. It is important to note that $A^*$
is a universal set, not depending on the initial state, in the sense that the `trap'  is $A^*\cap \{ k+1,\cdots\}$ for  proceeding from position $k$. 
This is different from the problem with  continuous time, where the optimal traps in different states $(t,k)$ are not consistent, unless  the point process of success epochs is Poisson.

\subsection{Tests for the monotone case of optimal stopping}

We proceed with the continuous time setting, assuming $\successprofile$ and $\priordist$ given.
In state $(t,k)$, the probability to isolate the last success with $z$-strategy is 
a convex mixture of Bernstein polynomials
\begin{eqnarray}\label{S11}
{\cal S}_1(t,k ;z):=
\sum_{j=1}^\infty 
\pi(j|t,k)
\sum_{i=0}^{j-1} {j\choose i}z^{i}(1-z)^{j-1}s_1(k+i+1,k+j).
\end{eqnarray}
The $z=0$ instance
$$
{\cal S}_1(t,k; 0)=\sum_{j=1}^\infty 
\pi(j|t,k) s_1(k+1,k+j),
$$
 is the probability to win with $\snext$ and ${\cal S}_1(t,k; 1)=0.$
Similarly, the probability that none of the successes is trapped by the $z$-strategy is
\begin{eqnarray*}
{\cal  S}_0(t,k;z ):=\sum_{j=0}^\infty \pi(j|t,k) 
\sum_{i=0}^{j-1} {j\choose i}z^{i}(1-z)^{j-1}s_0(k+i+1,k+j),
\end{eqnarray*}
and ${\cal  S}_0(t,k;0)$ is the probability to win with $\sbygone$.\jump

Using (\ref{poster}) and (\ref{x-t}), we can cast the winning probabilities as
\begin{eqnarray*}
{\cal S}_0(t,k;0)&=&f_k(x) P_k(x),\\
{\cal S}_1(t,k;0)&=&f_k(x) Q_k(x),\\
{\cal S}_1(t,k; z)&=&f_k(x) R_k(x,z),
\end{eqnarray*}
where
\begin{eqnarray*}
P_k(x)&:=&\sum_{j=0}^\infty {k+j\choose j}w_{k+j} x^j s_0(k+1,k+j),\\
Q_k(x)&:=&\sum_{j=1}^\infty {k+j\choose j}w_{k+j} x^j s_1(k+1,k+j),\\
R_k(x,z)&:=&\sum_{j=1}^\infty {k+j\choose j}w_{k+j} x^j  \sum_{i=0}^{j-1} {j\choose i} z^i (1-z)^{j-i}  s_1(k+i+1,k+j).
\end{eqnarray*}
Thus, $Q_k(x)=R_k(x,0)$. We are looking next at some critical `cutoffs'  for the trapping game and optimal stopping.

\begin{lemma}\label{L2} Equation $P_k(x)=Q_k(x)$ has at most one root $\alpha_k>0$, for every $k\geq 1$.
\end{lemma}
\begin{proof}
The series $P_k(x)-Q_k(x)$ has at most one change of sign from $+$ to   $-$, hence Descartes' rule of signs  for power series \cite{Curtiss} entails that there is at most one root.
\end{proof}
We set $\alpha_k=\infty$ if the root does not exist. Define the cutoff
$$
a_k=\left(1-\frac{1-\alpha_k}{q}\right)_+.
$$ 
This is the earliest time when $\sbygone$ becomes as beneficial as $\snext$.     
The {\it myopic} stopping strategy starting at time $t_0$
 is defined as
$$\tau^*:=\inf\{t\geq t_0:  (t,N_t) ~{\rm is~a~success~epoch~and~} t\geq a_{N_t}\}.$$ 
Keep in mind that if the sequence $(\alpha_k)$ is monotone, then $(a_k)$ is also monotone but with the monotonicity direction reversed.
The monotone case of optimal stopping and optimality of  the myopic strategy hold if $a_k\downarrow$.

\begin{lemma}\label{L3}  Equation $D_z R(x,0)=0$  has at most one root $\beta_k>0$, for every $k\geq 0$. If the root exists, then $\beta_k\leq\alpha_{k+1}$. 
\end{lemma}

\begin{proof} 
We follow the argument in Lemma \ref{L2}.
The derivative at $z=0$ is
$$  D_z R_k(x,0)=p _{k+1}\,\sum_{j=1}^\infty {k+j\choose j}w_{k+j}\,
j\,x^j\,\{s_0(k+2,k+j)-s_1(k+2,k+j)\}.$$
 This has at most one change of sign as $x\geq 0$ varies, and then from $+$ to $-$.
Furthermore, 
\begin{eqnarray*}
D_z R_k(x,0) &\geq&p_{k+1} \sum_{j=1}^\infty  {k+j\choose j} w_{k+j} x^j
\{s_0(k+2,k+j)-s_1(k+2,k+j)\}\\
&=&p_{k+1}\{ P_{k+1}(x) -Q_{k+1} (x)\}.
\end{eqnarray*}
This   follows by comparing the series and noting that 
 the weights at positive terms in $D_z$ are higher.
\end{proof}

If there is no finite root, we set $\beta_k=\infty$.
Let
$$
b_k:=\left(1-\frac{1-\beta_k}{q}\right)_+.
$$ 
We have $D_z R_k(x, 0)<0$
for  $ t\in(b_k,1]$, and $b_k\geq a_{k+1}$ by  Lemma \ref{L3}. Thus, $b_k$ is the earliest time when the action 
{\tt next}at index $k$ cannot be improved  by a $z$-strategy with small $z$.\jump

To summarise the above: 
\begin{itemize}
	\item For $t<a_k$: $\snext$ is better than $\sbygone$.
	\item For $t<b_k$: a trapping strategy is better than {\tt next}.
\end{itemize}

\begin{thm}\label{thma}  The optimal stopping problem belongs to the monotone case (for every $q$ and arbitrary initial state) if and only if
$\alpha_1\leq\alpha_2\leq\cdots.$ In that case we have the interlacing pattern of roots
\begin{equation}\label{interl}
\cdots\leq \alpha_k\leq \beta_k\leq\alpha_{k+1}\leq\beta_{k+1}\leq\cdots
\end{equation}
\end{thm} 
\begin{proof} 
We argue in probabilistic terms.
The bivariate sequence of success epochs $(t,k)^\circ$ is an increasing Markov chain. The monotone case of optimal stopping occurs iff  the set 
of states where $\sbygone$ outperforms $\snext$ is closed, which holds iff this is an upper subset 
with respect to the partial order 
in $[0,1]\times\{1,2,\cdots\}$. 
The latter property amounts to the monotonicity condition $\alpha_k\uparrow$.

By Lemma 3, the inequality $\alpha_k\leq \beta_{k+1}$ always hold.  In the monotone case, 
if in some state  $(t,k)^\circ$ the actions {\tt bygone} and {\tt next} are equally good, then trapping cannot improve upon these by optimality of the myopic strategy. 
In the analytic terms,  the above translates as the inequality $\beta_k\leq\alpha_k$.
\end{proof}

The monotone case does not hold if $\alpha_{k+1}<\alpha_k$ for some $k$.
Alternatively, one can use $\beta_k<\alpha_k$ as a test. Indeed, 
If  $\beta_k<\alpha_k$ and $a_k\neq0$ then 
in state $(a_k, k)^\circ$ a trapping strategy is better than both $\sbygone$ and $\snext$.

\subsection{Unimodality and concavity}\label{UC}

Being a convex mixture of unimodal functions,  ${\cal S}_1(t,k;\cdot)$ itself need not be unimodal.
Accordingly, the optimal trap may not be a final interval. 
Concavity is a simple condition to ensure the unimodality of ${\cal S}_1(t,k;\cdot)$.\jump

Suppose $s_1(\cdot,n)$ is concave for every $n\geq 1$, that is, the second difference in the first variable is non-positive.
By the shape-preserving properties of Bernstein polynomials, the internal sum  in (\ref{S11}) is a concave function of $z$, hence the mixture ${\cal S}_1(t,k;\cdot)$ is also concave. In that case we have

\begin{thm}\label{T2} If $s_1(\cdot,n)$ is concave for every $n$, then for cutoff
$z$  coinciding with  the mode of ${\cal S}_1(t,k;\cdot)$,
the $z$-strategy is optimal among all trapping strategies. 
The mode is distinct from $0$ 
iff $t<b_k$.
\end{thm}
\begin{proof} The overall 
optimality  follows from the unimodality as in Theorem \ref{T1}.
By concavity, the mode is zero if $D_z(t,k,\priordist, 0)\leq0$, and is positive otherwise.
\end{proof}

The concavity is easy to express in terms of $\successprofile$ explicitly. For instance, consider the second difference for $k=1$.
 The second difference in the variable $k$ of the probability generating function 
$$\lambda\mapsto \prod_{j=k}^n (1-p_j+\lambda p_j)$$
becomes
$$\{ (1-p_k+\lambda p_k)(1-p_{k+1}+\lambda p_{k+1})-2(1-p_{k+1}+\lambda p_{k+1})+1    \}\prod_{j=k+2}^n (1-p_j+\lambda p_j).$$
Computing $D_\lambda$ at $\lambda=0$ yields the second difference of $s_1(\cdot\,,\,n)$
\begin{equation}\label{pkk}
(p_k-2 p_k p_{k+1}-p_{k+1})+(p_kp_{k+1}-p_k+p_{k+1})\sum_{j={k+2}}^n \frac{p_j}{1-p_j}.
\end{equation}
Hence, a
 sufficient condition for the concavity of  $s_1(\cdot,n)$ is 
\begin{equation}
\label{concave}
p_k-2p_k p_{k+1}-p_{k+1}\leq 0, ~~~p_k p_{k+1}-p_{k}+p_{k+1}\leq 0, ~~~~k\geq 1.
\end{equation}

\begin{figure}
\centering{
\includegraphics[scale=0.8]{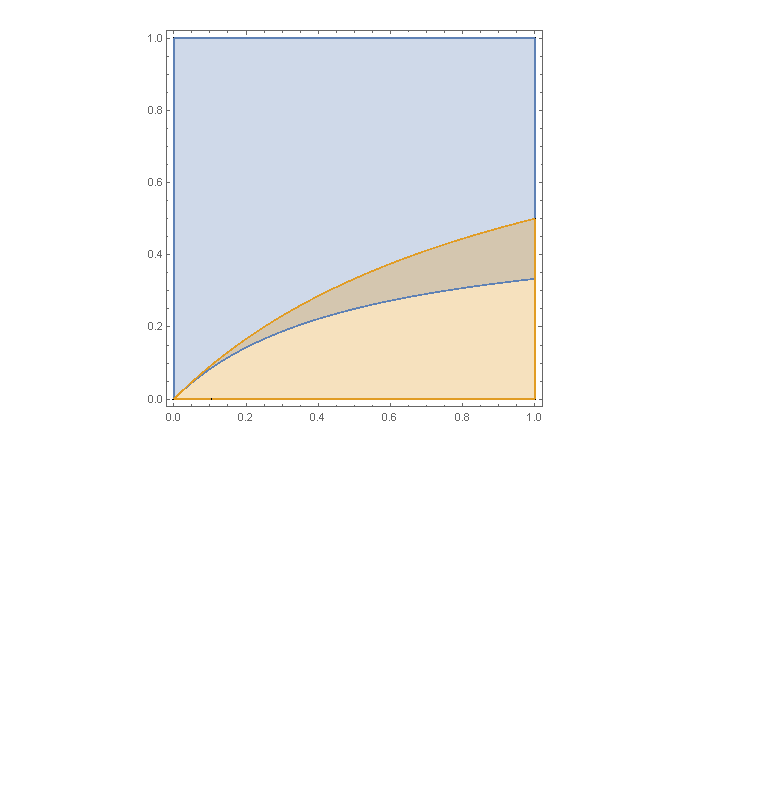}}
\vskip-5cm
\caption{Concavity  (\ref{concave}) holds for the region between parabolas.}
\label{fig:conacave}
\end{figure}

We stress that  (\ref{concave}) 
 ensures unimodality for arbitrary $\boldsymbol\pi$ and only involves two consequitive success probabilities.
The price to pay for the generality  is that the condition is  restrictive, as seen on Figure \ref{fig:conacave}.

For the {\it Karamata-Stirling} profile, straight calculation shows that (\ref{pkk}) is non-positive, hence  $s_1(\cdot,n)$'s concave,
 iff
$$ \frac{1}{2}\leq \theta\leq 1.$$
This is a narrow range, but it includes two most important cases $\theta=1$ and $\theta=1/2$.

\section{The best choice problem under the  log-series prior}

In this section, we consider the classic profile $p_k=1/k$ from the random records model
and  the logarithmic series prior 
\begin{equation}\label{logser}
\pi_n=  c(q)\,\frac{q^n}{n}, ~~~n\geq 1,
\end{equation}
(so $\pi_0=0$), where
$c(q)=|\log(1-q)|^{-1}$.
Two representations  of such $\boldsymbol\pi$ as mixed Poisson distribution can be obtained by mixing
 before  zero-truncating or after \cite{Johnson}.

From a wider view, the setting  is the $(\nu=0, \theta=1)$  instance of the  problem with negative binomial prior NB$(\nu,q)$ and the Karamata-Stirling profile as studied in \cite{GD}.
Similarly to the case  $0<\nu<1$, comparison with the geometric prior yields 
$$\alpha_k\geq 1-1/e, ~~~\alpha_k\to 1-1/e.$$
This entails that the  roots sequence $(\alpha_k)$ cannot be decreasing, hence by Theorem \ref{thma}  the monotone case of optimal stopping does not hold.
As was shown in Section \ref{UC} for this profile, the best trapping strategy is a $z$-strategy.

Let $T_1$ be the time of the first trial.

\begin{lemma}\label{log-ser-set} Under the logarithmic series prior {\rm\, (\ref{logser})}
\begin{itemize}
\item[\rm(i)] The time of the first trial  $T_1$ has probability density function
$$ t\mapsto  \frac{c(q)\,q}{1-(1-t)q},~~~t\in[0,1].$$
\item[\rm(ii)]
$(N_t,~t\in[0,1])$ is a P{\'o}lya-Lundberg birth process with transition rates
\begin{eqnarray*}
\prob(N_{t+{\rm d}t}-N_t=1\,| \,N_t=k)= 
\begin{dcases}
\frac{c((1-t)q)\,q}{1-(1-t)q},~~~k=0,\\
\frac{k}{t+q^{-1}-1},~~~k\geq 1,
\end{dcases}
\end{eqnarray*}
\item[\rm(iii)] Given  $N_t=k$, the posterior distribution $\pi(\cdot\,|\,t,k)$ of $N_1-N_t$ is 
${\rm NB}(k, (1-t)q).$
In particular, conditionally on $T_1=t_1$, the posterior distribution is geometric with  the `failure' probability $(1-t_1)q$.

\end{itemize}
\end{lemma}
\begin{proof}
Assertion (i) follows from
$${\mathbb P}(T_1>t)= {\mathbb P}(N_t=0)=
 \sum_{n=1}^\infty \frac{c(q) q^n (1-t)^n}{n\,},$$
and
 (iii) from the identity
$${k+j\choose j}\frac{x^j}{k+j}={k+j-1\choose j}\frac{x^j}{k}$$
underlying $\pi(j|t,k)$ for $x=(1-t)q$.
\end{proof}

The value $q=1$ is on the edge of convergence. It formally corresponds to the infinite `non-informative' prior. 
As a result of that, the P{\'o}lya-Lundberg process is well defined by the rates in  (ii) for any initial state $(t_0,k_0)$ with  $t_0>0$.
With initial state $(t_0,0)$, the model is equivalent to the model with logarithmic prior {\rm NB}$(0,(1-t_0))$ and trials occurring on $[0,1]$.
In the  $t_0\to 0$ limit,  the process of record times becomes a  Poisson process with intensity function $t^{-1}$ with the $1/e$-strategy being then  optimal.

\subsection{Hypergeometrics}

It will be helpful to recall some properties 
 of the Gaussian hypergeometric function
 $$
F(a,b,c;x):=\sum_{j=0}^\infty \frac{(a)_j(b)_j}{(c)_j}\frac{x^j}{j!}
$$
These include: the differentiation formula
$$D_x F(a,b,c,x)=\frac{ab}{c}\,F(a+1,b+1,c+1,x),$$
the transformation formula
$$F(a,b,c;x)=(1-x)^{c-a-b}F(c-a,c-b,c;x),$$
and Euler's integral representation for $c>b>0$
$$F(a,b,c;x)=\frac{\Gamma(c)}{\Gamma(b)\Gamma(c-b)}\int_0^1 \frac{y^{b-1}(1-y)^{c-b-1}{\rm d}y}{(1-xy)^a}.$$

The probability  generating function for the number of successes following state  $(t,k)$, for $k\geq 1$, is  given by a hypergeometric function
\begin{eqnarray*}
\lambda\mapsto   (1-x)^k  \sum_{j=0}^\infty {k+j-1\choose j} x^j \frac{(k+\lambda)_j}{(k+1)_j}&=&\\
(1-x)^k \sum_{j=0}^\infty \frac{(k)_j(k+\lambda)_j}{(k+1)_j}\frac{x^j}{j!}&=&\\
(1-x)^k \,F(k+\lambda, k, k+1;x).
\end{eqnarray*}
We read off that the normalisation function is $f_k(x)=k(1-x)^k$ for $k\geq1$, and $f_0(x)= |\log(1-x)|^{-1}$.
Expanding at $\lambda=0$ we identify two basic power series as
\begin{eqnarray*}
P_k(x)&=&  k^{-1}\, F(k, k, k+1; x),\\
Q_k(x)&=& k^{-1}\,D_a F(k,k,k+1;x),
\end{eqnarray*} 
where  as before   $x=(1-t)q\in[0,1]$ and
$D_a$ is the derivative in the first parameter. 
The differentiation formula implies backward recursions
\begin{eqnarray}
\nonumber
D_x P_k(x)&=&k P_{k+1}(x),\\ 
\label{recQ}
D_x Q_k(x)&=&P_{k+1}(x)+k\, Q_{k+1}(x).
\end{eqnarray}

Applying  the transformation formula yields
$P_k(x)=(1-x)^{1-k} F(1,1,k+1,x),$ 
hence, we may write the winning probability with $\sbygone$ as the series
\begin{eqnarray*}
{\cal S}_0(t,k;0)
=(1-x)\sum_{j=0}^\infty  \frac{j!\, x^j} {(k+1)_j},~~~~x=(1-t)q.
\end{eqnarray*}
It is readily seen that as $k$ increases, this function decreases to $1-x$.
The fact was shown in 
  \cite{BrussRogers}  probabilistically.
Convergence to $1-x$ is related to the fact that for large $k$ the process of record times approaches  a Poisson process.

Explicitly, for $k=1, 2$ and  $L:= -\log(1-x)$,   computing the roots to six decimal places we have
\begin{eqnarray*}
P_1(x)=  \frac{L}{x}\,, ~~~~~~~~~~~~~~~~~Q_1(x)=\frac{ L^2}{2 x},~~~~~~~~~~~~~~~~~~~~~~~~~~~~~\alpha_1&=&1-e^{-2}=0.864665,\\
P_2(x)=\frac{2( x-L +x L)}    {(1-x)x^2},~
Q_2(x)=     \frac{-2 x+2L -L^2+xL^2}{(1-x)x^2}, ~~~~\alpha_2&=&0.755984\jump
\end{eqnarray*}
\begin{assertion} The roots satisfy $\alpha_k\downarrow \,1-1/e$ as $k\to\infty$. 
\end{assertion}
\begin{proof}
In the case of constant weights $w_j=1$, the prior is geometric and  all roots coincide with $1-1/e$.
The log-series distribution weights satisfy $w_{n+1}/w_n\uparrow 1$, hence, comparison with the geometric distribution (see \cite{GD}) gives $\alpha_k>1-1/e$ and 
$\alpha_k\to 1-1/e$. That the sequence of roots is decreasing will be shown separately.
\end{proof}

\begin{cor}\label{cor8} The optimal stopping problem is not monotone, the myopic strategy $\tau^*$ is not optimal, and
\begin{itemize}
\item[\rm(i)] for  $q> 1-1/e$ the myopic strategy  is determined by an infinite sequence of cutoffs converging to $1-(1-1/e)/q$.
\item[\rm(ii)] for $t\geq (1-(1-1/e)/q)_+$,  $\tt bygone$ is the optimal action  for every $k$,
\item[\rm(iii)] 
for times as in {\rm (ii)}
the optimal stopping strategy stops greedily at the first available record.
\end{itemize}
\end{cor}

With some manipulation, we can derive an integral formula for $R(x,z)$.
Consider first $k\geq 1$.
The probability generating function of the number of record epochs following $(t,k)$ and falling
in the final interval $[t+z(1-t),1]$ has probability generating function
\begin{eqnarray*}
\lambda\mapsto (1-x)^k\sum_{j=0}^\infty {k+j-1\choose j}x^j \sum_{i=0}^j {j\choose i} z^i(1-z)^{j-i} \frac{(k+i+\lambda)_{j-i}}{(k+i+1)_{j-i}}& =&\\
(1-x)^k\sum_{i=0}^\infty {k+i-1\choose i} (xz)^i F(k+i+\lambda, k+i, k+i+1;x-xz)&=&\\
k(1-x)^k\sum_{i=0}^\infty {k+i\choose i} (xz)^i  \int_0^1 \frac{y^{k+i-1} {\rm d}y}{(1-xy+xyz)^{k+i+\lambda}}&=&\\
k(1-x)^k\,\int_0^1 \frac{y^{k-1} (1-xy+xyz)^{1-\lambda}{\rm d}y}{(1-xy)^{k+1}}.
\end{eqnarray*}
Differentiating at $\lambda=0$
yields for $x=(1-t)q,~z\in [0,1]$
\begin{equation}\label{Rk}
{\cal S}_1(k,t;z)=k(1-x)^k R_k(x,z)=
k(1-x)^k\int_0^1 \frac{y^{k-1}   (1-xy+xyz)   | \log (1-xy+xyz) |                   {\rm d}y}      {(1-xy)^{k+1}}.
\end{equation}
For $k=0$, a similar calculation  with log-series weights NB$(0, x)$ gives for $x=(1-t)q$
$$
{\cal S}_0(0,t,z)= \frac{ R_0(x,z) }{|\log(1-x)|}  = \frac{1}{\log(1-x)}    \int_0^1 \frac{(1-x y+xyz)\log(1-xy+xyz)}{y(1-x y)}\, {\rm d}y.
$$

\subsection{Monotonicity of cutoffs for the myopic strategy}

We show next that the roots  are indeed decreasing, which is the direction opposite to the one needed for optimality of the myopic strategy.
We may define the root $\alpha_k$ in terms of the quotient, as a 
 unique solution on $[0,1)$  to
\begin{equation}\label{quot}
\frac{Q_k(x)}{P_k(x)}=1~~~\Longleftrightarrow~~~\frac{D_a F(k,k,k+1;x)}{F(k,k,k+1;x)} =1.
\end{equation}
As $x$ runs from $0$ to $1$, the quotient varies from $0$ to $\infty$.

Euler's integral  for the hypergeometric function  specialises as
$$
F(k+\lambda,k,  k+1; x)=k \int_0^1 \frac{y^{k-1}{\rm d}y}{(1-xy)^{k+\lambda}}{\rm d}y.
$$
Expanding at $\lambda=0$ gives
\begin{equation}\label{EulInt}
P_k(x)= k \int_0^1 \frac{y^{k-1}   {\rm d}y}{(1-xy)^{k}},~~~~~~~~
Q_k(x)= k \int_0^1 \frac{y^{k-1} |\log(1-x y)|    {\rm d}y}{(1-xy)^{k}}.
\end{equation}

\begin{lemma} The logarithmic derivative {\rm(\ref{quot})} increases in $k$,  hence $\alpha_k \downarrow (1-1/e)$.
\end{lemma}

\begin{proof} From (\ref{EulInt})
\begin{eqnarray*}
Q_k(x)P_{k+1}(x)=  \int_0^1 \frac{y^{k-1} |\log(1-x y)|    {\rm d}y}     {(1-xy)^{k}} 
\int_0^1 \frac{z^{k}  {\rm d}z}     {(1-xz)^{k+1}}&=&\\
\int_0^1\int_0^1 \frac{y^{k-1}z^{k-1} |\log(1-x y)|  } {(1-xy)^{k}(1-xz)^{k}} 
 \frac{z}   {(1-xz)}       {\rm d}y  {\rm d}z,
\end{eqnarray*}
By the same argument, a  similar formula is obtained for $Q_{k+1}(x)P_{k}(x)$. Splitting the integration domain,
then swapping the variables on the triangle above the diagonal $y>z$ yields
\begin{eqnarray*}
{Q_k(x)}    {P_{k+1}(x)}- {Q_{k+1}(x)}{P_{k}(x)}&=&\\
\int_0^1\int_0^1 \frac{y^{k-1}z^{k-1} |\log(1-x y)|  } {(1-xy)^{k+1}(1-xz)^{k+1}} 
\left( z-y    \right) {\rm d}y  {\rm d}z&=&\\
\int\!\!\!\!\!\!\!\!\! \int\limits_{0<y<z<1}
 \frac{y^{k-1}z^{k-1} } {(1-xy)^{k+1}(1-xz)^{k+1}} 
 \log\left(\frac{1-x z}{1-xy }\right)
\left( z-y    \right) {\rm d}y  {\rm d}z&<&0,
\end{eqnarray*}
because the symmetric part of the integrand is positive and the asymmetric is negative for  $x\in[0,1)$.
\end{proof}

\subsection{The information bounds}

Suppose that in state $(t,k)$ the gambler learns that there are {\it exactly $j$ trials yet to occur}. A higher winning probability is attainable with more information and it is one of:

\begin{itemize}
\item[(i)] $s_0(k+1,k+j)$ for {\tt bygone},
\item[(ii)] $s_1(k+1, k+j)$ for {\tt next},
\item[(iii)] $\max_{z\in[0,1]}  S_1(k, k+j,z)$ for the best trapping,
\item[(iv)] $\max_{k': k'>k} s_1(k'+1,k+j)$ for the optimal stopping strategy, now independent of the time of trials.
\end{itemize}

Weighting these with the posterior distribution $\pi(j\,|\,t,k)$ gives upper bounds $I_k$, on the winning probability, only achievable by the informed gambler.


\begin{figure}[h]
	\centering
	\subfloat{{\includegraphics[scale=0.28]{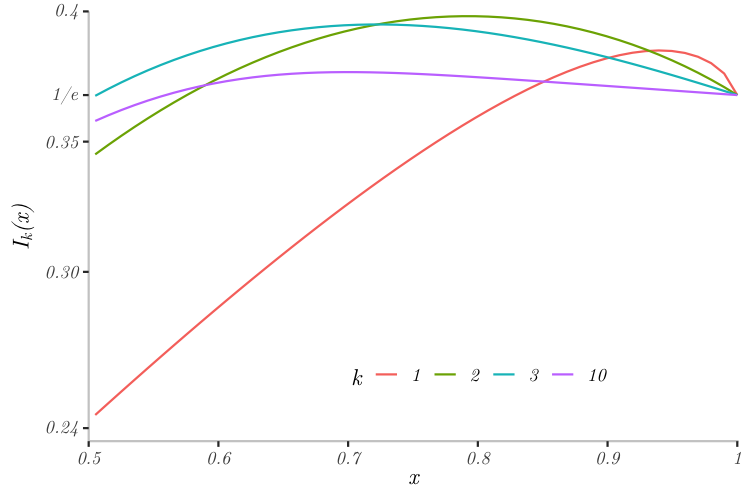} }}%
	\,
	\subfloat{{\includegraphics[scale=0.28]{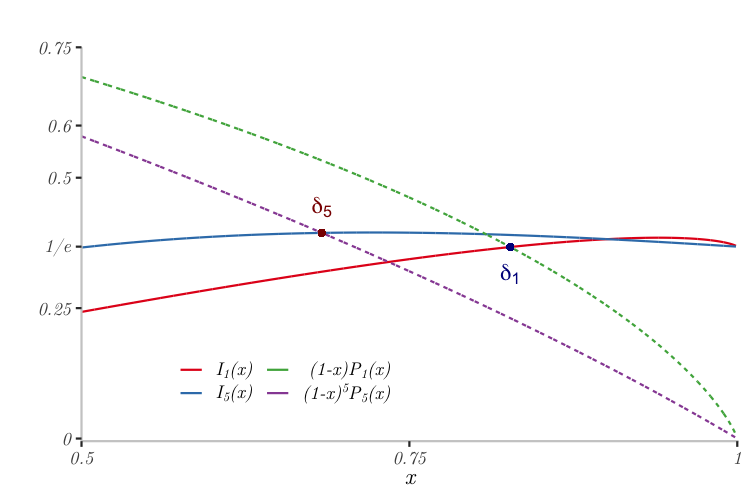} }}%
	\caption{Bounds on the optimal strategy $I_k(x)$}
	\label{fig:infobound}%
\end{figure}

\subsection{The value function}

Define $v(t,k)$ to be the {\it continuation value} of state $(t,k)$, equal to the winning probability achieved by the  optimal stopping strategy    
starting in this state. By the optimality principle, in state $(t,k)^\circ$ it is optimal to stop (action {\tt bygon}) iff ${\cal S}_0(t,k;0)\geq v(t,k)$.
We have $v(1,k)=0$ for $k\geq 1$ and $v(1,0)=1$ since, near the end of the time interval, it is unlikely to see more trials if some have occurred,
but at least one trial is ensured by the log-series prior if none occurred. Passing to $x=(1-t)q$ we can write the continuation value as a function 
$V_k(x)$ 
so that 
$$v(t,k)=V_k((1-t)q), ~~~k\geq 0.$$

The optimality principle yields a recursion for the $V_k$'s as follows:\jmp{3}
Given $N_t=k$, let $T_{k+1}$ be the next trial epoch or $1$ in the event $N_1=k$.
Similarly to the argument in Lemma \ref{log-ser-set}, it is seen that the random variable
$(1-T_{k+1})/(1-t)$
has density:
$$y\mapsto \frac{kx(1-x)^k}{(1-x+xy)^{k+1}},~y\in(0,1].$$
By the \nth{(k+1)}{st} trial, the optimal stopping strategy chooses a better action in case the trial is a success, hence integrating out $T_{k+1}$ we obtain
$$
V_k(x)=\int_0^1\left[ \frac{1}{k+1}\max \{(1-y)^{k+1}P_{k+1}(y), V_{k+1}(y)\}+\frac{k}{k+1} V_{k+1}(y) \right] \frac{kx(1-x)^k{\rm d}y}{(1-x+xy)^{k+1}}.
$$
 This has the following differential form for $k\geq 1$
\begin{equation}\label{DF1}
(1-x)\,D_x V_k(x)=\frac{k}{k+1} \left((1-x)^{k+1}P_{k+1}(x)-V_{k+1}(x)  \right)_+ +k\{V_{k+1}(x)-V_k(x))\}
\end{equation}
The instance $k=0$ is special. 
Integrating out the variable $T_1$ gives:
$$V_0(x)=\int_0^1 \max((1-y)P_1(y),V_1(y)) \frac{{\rm d}y}{(1-x+xy)|\log(1-x)|},$$  
or, in the differential form with initial conditions $V_0(0)=1$ and  $V_k(0)=0 \text{ for} ~k\geq 1$:
\begin{equation}\label{DF2}
(1-x)|\log(1-x)|\, D_x V_0(x)= \max\{(1-x)P_1(x), V_1(x) \}- V_0(x),
\end{equation}

By the theory of optimal stopping  \cite{CRS}, 
the value   function can be characterised as the minimal solution to  (\ref{DF1}), (\ref{DF2}). 
For computational purposes, one can use 
 the limit relation
$$\lim_{k\to\infty} V_k(x)= \max\{(1-x)| \log(1-x)|, ~1/e\}$$
as a boundary condition at  $k=\infty$.\jump

On the left part of the interval we know the value function exactly
\begin{eqnarray}
V_k(x)= k(1-x)^k Q_k(x), ~~~{\rm for~}0\leq x\leq 1-1/e, ~~k\geq 0,
\end{eqnarray}
as a consequence of Corollary \ref{cor8}. As a check, for $k\geq 1$ let $\widehat{V}_k(x):=k^{-1}(1-x)^{-k}V_k(x)$. With this change of variable, (\ref{DF1}) simplifies to
$$D_x \widehat{V}_k(x)=(P_{k+1}(x)-\widehat{V}_{k+1}(x))_+     + (k+1) \,\widehat{V}_{k+1}(x),$$  
For $x$ in the range where $P_{k+1}(x)-\widehat{V}_{k+1}(x)\geq 0$, this becomes the recursion (\ref{recQ}).

\begin{figure}[h]
	\centering
	\subfloat{{\includegraphics[scale=0.27]{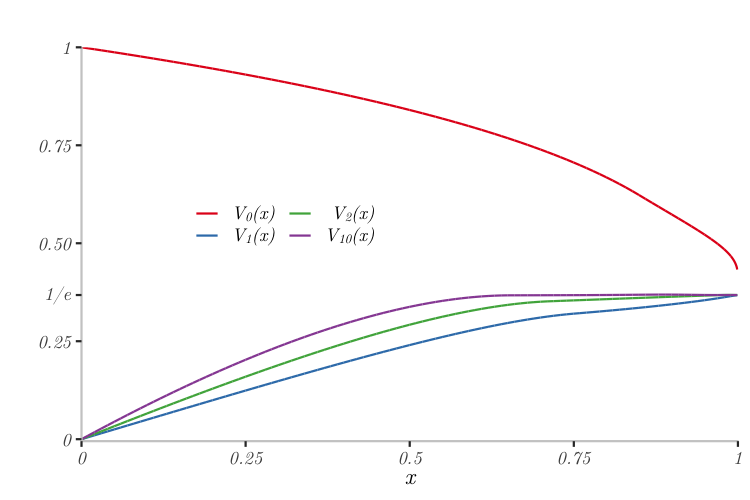} }}%
	\,
	\subfloat{{\includegraphics[scale=0.27]{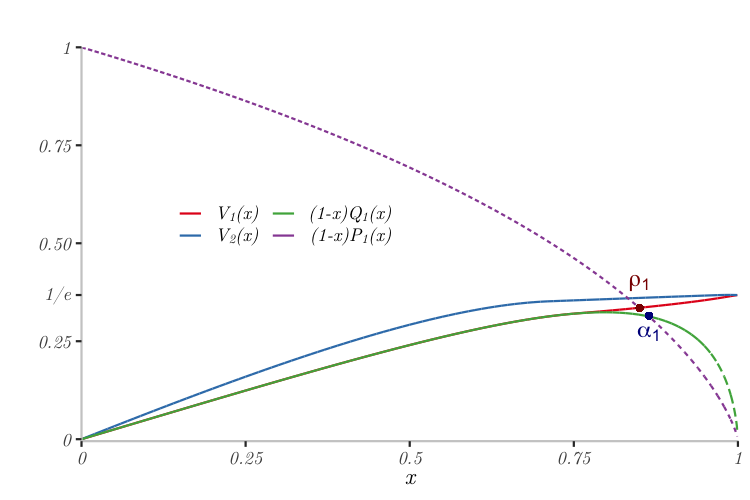} }}%
	\caption{Stop and Continuation Values}
	\label{fig:stopcont}%
	\label{F4}		
\end{figure}

To solve (\ref{DF1}), (\ref{DF2}) numerically for $1-1/e\leq x\leq$, one can use the endpoint values
$$V_k(1-1/e)= k e^{-k} Q_k(1-1/e)$$ 
along with the convergence
$$\lim_{k\to\infty} V_k(x)=\begin{cases}~~~ 1/e,~~~~~~~~~~~~~{\rm for~}1-1/e\leq x\leq 1,\\
-(1-x) \log(1-x), ~~~{\rm for~}0\leq x\leq 1-1/e,
\end{cases}
$$
in the role of a boundary condition at  $k=\infty$.
Figure \ref{F4} shows some shapes.

Define:
$$\gamma_k:=\sup\{x: V_k(x)\geq (1-x)^k P_k(x)\} , ~k\geq 1.$$
Since $\gamma_k<\alpha_k$, we have $\gamma_k\to 1-1/e$. 
It is natural to expect that $\gamma_k$'s are decreasing, and that the optimal stopping strategy is determined by the cutoffs $(1-\gamma_k/q)_+$, in complete analogy with the myopic strategy $\tau^*$.

This is confirmed by  simulation which also shows that the myopic and optimal strategies are very close to one another, as is evident by comparing the critical  points in Table \ref{tab:roots}.
We remind that these are related to real-time cutoffs via (\ref{x-t}).

While $V_k(x)\to 1/e$  when either $x\to 1$ or $k\to \infty$, the simulation shows that the functions increase in $x$ and $k\geq 1$.
In contrast to the above, $V_0$ is decreasing with $V_0(x)\downarrow 1/e$ as $x\uparrow 1$, see Figure \ref{fig:stopcont}.
The latter convergence is slow, because the logarithmic distribution of $N_1|N_t=0$   puts a relatively  high weight on small values of $n$, which is advantageous for 
stopping at the last record.
For instance, for $q=1-10^{-6}$ the mean is about $72382$ while the  probability of only one trial is still higher than  $0.072382$.

	\begin{table}[h]
		\centering

			\begin{tabular}{ |c|c|c|c|c|c| }
				\hline
				$k$ &$\alpha_k$& $\beta_k$& $\gamma_k$ & $\delta_k\tnote{2}$ & $\rho_k$\tnote{3} \\
				\hline
				1 & 0.864665&  0.756004 & 0.849635 &  0.826893& 0.850335  \\
				2 & 0.755984&  0.714616 & 0.753621  &  0.718332   & 0.753727 \\
				3 & 0.714596&  0.693549& 0.713957 & 0.683295   & 0.713995\\       
				4 & 0.693529&  0.680931 & 0.693375&0.668986&  0.693311 \\
				5 & 0.680911 &0.672567 & 0.680887 & 0.661520&0.680814\\
				10 & 0.656034 & 0.653833 & 0.656109& 0.647653& 0.656028\\
				\hline

			\end{tabular}
			
	\vskip0.5cm
\captionsetup{singlelinecheck=off, format=plain, font=small, justification=justified}

\caption[foo]{
{\centering
	Critical points:

		\begin{itemize}[leftmargin=5cm, font=\small\color{blue},itemsep=0pt, parsep=0pt,]
		\item[{\rm $\alpha_k$:}] 
		critical points for the myopic strategy
		\item[{\rm $\beta_k$:}] 
		balance points where $\snext$ is  as good as trapping
		\item[{\rm $\gamma_k$:}] 
		critical points for the optimal strategy
		\item[{\rm $\delta_k$:}] 
		lower bounds of $\gamma_k$ obtained from the information bound
		\item[{\rm $\rho_k$:}] 
		balance points where $\sbygone$ is as good as trapping
	\end{itemize}

}}

		\label{tab:roots}
	\end{table}

\end{document}